    \documentclass[12pt]{article}
    \usepackage{t1enc}
    \usepackage[latin1]{inputenc}
    \usepackage[english]{babel}
        \usepackage{latexsym}
        \newcommand{\bp}{{\bf P}}

    \newtheorem{conjecture}{Conjecture}[section]
    
\newcommand{\be}{{\bf E}}
\newtheorem{theorem}{Theorem}[section]
\newtheorem{corollary}{Corollary}[section]
\newtheorem{lemma}{Lemma}[section]

\def\begg{\begin{equation}}
\def\endd{\end{equation}}
\def\pe{{\bf P}}
\def\ep{{\varepsilon}}
\begin{document}

\centerline{\Large\bf ON THE NUMBER OF CUTPOINTS OF}

\medskip
\centerline{\Large\bf THE TRANSIENT NEAREST NEIGHBOR}

\medskip
\centerline{\Large\bf  RANDOM WALK ON THE LINE}

\bigskip \bigskip \bigskip \bigskip \bigskip

\renewcommand{\thefootnote}{1} \noindent
\textbf{Endre Cs\'{a}ki}\footnote{Research supported by the
Hungarian National Foundation for Scientif\/ic Research, Grant No.
  K 61052 and K 67961.}\newline Alfr\'ed R\'enyi Institute of
Mathematics, Hungarian Academy of Sciences, Budapest, P.O.B. 127,
H-1364, Hungary. E-mail address: csaki@renyi.hu

\bigskip

\renewcommand{\thefootnote}{2} \noindent \textbf{Ant\'{o}nia
F\"{o}ldes}\footnote{Research supported by a PSC CUNY Grant, No.
61061-0039.}\newline Department of Mathematics, College of Staten
Island, CUNY, 2800 Victory Blvd., Staten Island, New York 10314,
U.S.A. E-mail address: foldes@mail.csi.cuny.edu

\bigskip

\noindent \textbf{P\'al R\'ev\'esz}$^1$ \newline Institut f\"ur
Statistik und Wahrscheinlichkeitstheorie, Technische Universit\"at
Wien, Wiedner Hauptstrasse 8-10/107 A-1040 Vienna, Austria. E-mail
address: reveszp@renyi.hu

\bigskip \bigskip \bigskip
\noindent
{\bf Abstract} We consider transient nearest neighbor random walks
on the positive part of the real line. We give criteria for the
finiteness of the number of cutpoints and strong cutpoints. Examples and
open problems are presented.

\bigskip


\bigskip

 \vspace{.1cm}

\noindent Running head: Cutpoints of NN random walk

\vfill
\renewcommand{\thesection}{\arabic{section}.}

\section{Introduction}

\renewcommand{\thesection}{\arabic{section}} \setcounter{equation}{0}
\setcounter{theorem}{0} \setcounter{lemma}{0}
Let $X_0=0,\ X_1,X_2,\ldots$ be a Markov chain with
\begin{eqnarray}\label{def}
E_i&:=&\bp(X_{n+1}=i+1\mid X_n=i)=1-\bp(X_{n+1}=i-1\mid X_n=i)\\
&=&\left\{\begin{array}{ll} 1\quad & {\rm if}\quad  i=0\\
\nonumber
 1/2+p_i\quad & {\rm if}\quad i=1,2,\ldots,
\end{array}\right.
\end{eqnarray}
where $-1/2< p_i< 1/2,\ i=1,2,\ldots$.

\medskip\noindent
 {\bf Theorem A} (\cite{CH}, page 74) {\it Let $X_n$ be a Markov chain
 with transition probabilities given in {\rm (\ref{def}).}
Define
\begg U_i:={\frac{1-E_i}{E_i}}={\frac{1/2-p_i}{1/2+p_i}},\qquad
i=1,2,\ldots.
\label{uif}
\endd
Then $X_n$  is transient if and only if}
\begg \sum_{k=1}^{\infty}
\prod_{i=1}^k U_i < \infty. \label{trans}\endd
In case $p_i\geq 0$
the sequence $\{X_i\}$ describes the motion of a particle which
starts at zero, moves over the nonnegative integers and going away
from 0 with a larger probability than to the direction of 0.  We
suppose throughout this paper that $0\leq p_i<1/2,\, i=1,2,\ldots.$

In \cite{CSFR} we introduced the quantities

\begin{eqnarray}\label{dem}
D(m,n):=\left\{\begin{array} {ll}
0 & {\rm if\ } n=m,\\
1 & {\rm if\ } n=m+1,\\
1+\displaystyle{\sum_{j=1}^{n-m-1}\prod_{i=1}^j U_{m+i}} & {\rm if\
} n\geq m+2
\end{array}\right.
\end{eqnarray}
and
\begg
\lim_{n\to \infty}D(m,n)=:D(m).
\label{dn}
\endd
Clearly (\ref{trans}) implies that if the walk is transient then $D(m)$
is finite for all $m=1,2,\ldots$.

The properties of this Markov chain, often called birth and death chain
were extensively studied. Some of these results are mentioned
e.g. in \cite{CSFR}. Our main concern in that paper was to study the
local time of $\{X_n\}$, defined by

\begg\xi(x,n):=\#\{k:0\leq k\leq n,\ X_k=x\},\quad x=0,1,2,\ldots,
\quad n=1,2,\ldots
\label{locn}\endd and \begg\xi(x):=\lim_{n\to
\infty}\xi(x,n). \label{loc1}
\endd

The first topic in that paper was to find upper class results for the
local time.

\medskip\noindent
{\bf Theorem B} {\it Assume that $p_R\to 0$ as $R\to\infty$. Then with
probability $1$ we have

\begg
  \xi(R)\leq 2(1+\varepsilon)
D(R)\log R
\label{loc2}\endd
for any $\ep>0$ if $R$ is large enough.

Moreover,
 \begg
 \xi(R)\geq M D(R) \quad  { \rm i.o.\,\,a.s.}
{\label{loc3}}\endd for any $M>0.$}

The next question was how small can the local time be. In particular, we
studied the number of sites $R$ where $\xi(R)=1$. We found that the
answer heavily depends on the sequence $\{p_R\}_{R=1}^\infty$.

We will say that the nearest neighbor (NN) random walk $X_n^*$ is {\it
slower} than $X_n$ (or equivalently, $X_n$ is {\it quicker} than
$X_n^*$) if
\begg
p_R^*\leq p_R\quad{\rm for\,\, all}\,\, R=1,2,\ldots
\label{slow}
\endd

It is obvious that the quicker is $X_n$, the more sites with local
time equal to 1 will occur.

\noindent
{\bf Remark 1.} In (\ref{slow}) the required inequality could be
relaxed to hold for all but finitely many $R$ only, since finitely
many $p_R$ have no effect on the asymptotic behavior of the walk.
The same remark applies throughout the paper, when we require certain
properties of the $\{p_R\}$ system.

Introduce the following notations:
\begin{eqnarray*}
\Lambda(1,i,B)&=&{\frac{B}{i}},\\
\Lambda(2,i,B)&=&{\frac{1}{i}}+{\frac{B}{i\log i}},\\
\ldots,\\
\Lambda(K,i,B)&=&{\frac{1}{i}}+{\frac{1}{i\log i}}+\ldots+{\frac{B}
{i\log i\log\log i\ldots\log_{K-1}i}},
\end{eqnarray*}
where $\log_r i=\log(\log_{r-1}i)$.

As a consequence of Theorem A, in \cite{CSFR} we obtained the following.
(See also a more general result given in Theorem 3 of Menshikov et al.
\cite{MAI}.)

\noindent
{\bf Fact 1} {\it If  for any $K=1,2,\ldots$

$$\displaystyle {p_i=\frac{\Lambda(K,i,B)}{4}},$$
then the Markov chain $\{X_n\}$ is recurrent if $B\leq 1$ and transient
if $B>1.$}

In the spirit of  Remark 1 above, it is enough if $p_i$ takes the
value given above with finitely many exceptions, but assuming that
$0\leq p_i<1/2$ for all $i=1,2,\ldots$.

We proved in \cite{CSFR} that if $\displaystyle
{p_i=\frac{\Lambda(1,i,B)}{4}}\,$ with
$B>1$, then we not only have infinitely many sites with local time 1,
but we have infinitely many increasing runs of sites each having local
time 1. More precisely we have

\noindent{\bf Theorem C}
 {\it Let $\{X_n\}$ be an NN random walk with
$\displaystyle{p_R=\frac{\Lambda(1,R,B)}{4}}=\frac{B}{4R}$ and
$ B>1$. Then with probability $1$ there exist infinitely many $R$ for
which
$$\xi(R+j)=1$$
for each $\displaystyle{j=0,1,2,\ldots,\left[\frac{\log\log R}{\log
2}\right].}$}

\bigskip
However, if $X_n^*$ is  transient but slower than $X_n$ in Theorem
C, then one might ask whether it still has infinitely many sites
with local time 1. It turns out that this is not always true. James
et al. \cite{JLP} proved a surprising result which implies the
following

\noindent {\bf Theorem D} {\it If $\{X_n\}$ is an NN random walk
with $\displaystyle{p_R=\frac{\Lambda(2,R,B)}{4}}$ and $B>1$, then
with probability $1$ $X_n$ has only finitely many sites $R$ with
$\xi(R)=1$.}

In fact, they formulated their results in terms of cutpoints.
Call the site $R$ a {\it cutpoint}\, if for some $k$, we have $X_k=R$
and $\{X_0,X_1 \dots X_k\}$ is disjoint from $\{X_{k+1},X_{k+2} \dots \},$
i.e. $X_i\leq R,\, i=0,1,\ldots,k$, $X_k=R$ and $X_i>R,\,
i=k+1,k+2,\ldots$

The original version of Theorem C in \cite{JLP} reads as follows.

\noindent
{\bf Theorem D$^*$}
{\it If $\{X_n\}$ is an NN random walk with
$$
\frac{c_1}{k(\log k)^{\beta}}\leq U_1U_2\ldots U_k\leq
\frac{c_2}{k(\log k)^{\beta}}
$$
for some $\beta>1$ and positive constants $c_1,\,c_2$, then $\{X_n\}$ is
transient and has only finitely many cutpoints a.s.}

Cutpoints and related intersection problems for more general stochastic
processes have been investigated extensively in the literature, starting
with Dvoretzky et al. \cite{DEK}, Erd\H os and Taylor \cite{ET60}. A nice
summary of this topic is given by Lawler \cite{LA}.

For usual random walk (sums of i.i.d. random variables) we mention the
following general result of James and Peres \cite{JP}, where the
definition of cutpoint is somewhat different from above. $S_k$ is called
a cutpoint there if $p(S_i,S_j)=0$ for all $(i,j)$ such that $0\leq
i<k<j$, where $p(x,y)$ is the one-step transition probability from $x$
to $y$.

\medskip\noindent
{\bf Theorem E} {\it Any transient random walk $\{S_k\}$ with bounded
increments on the lattice $Z^d$ has infinitely many cutpoints a.s.}
\medskip

To formulate our main result, we introduce the following definitions.

Call the site $R$ a {\it strong cutpoint}\, if for some $k$, we have
$X_k=R$, $X_i<R,\, i=0,1,\ldots, k-1$ and $X_i>R,\,
i=k+1,k+2,\ldots$. Observe that $R$ is a strong cutpoint if and only
if $\xi(R)=1$, i.e. if it is visited exactly once. Clearly every strong
cutpoint is a cutpoint, but not conversely. In words, $R$ is a cutpoint
if upcrossings from $R$ to $R+1$ occurs only once, but $R$ can be
visited from $R-1$ several times.

In this paper we give a criteria for a transient NN random walk
which determines whether the number of cutpoints (or strong
cutpoints) is finite or infinite almost surely.

\begin{theorem}
Let $X_0=0,X_1,X_2,\ldots$ be a transient Markov chain with transition
probability $E_i$ as in {\rm (\ref{def})} and $0\leq p_i< 1/2,\,
i=1,2,\ldots$ Let $D(n),\, n=1,2,\ldots$ be as in {\rm (\ref{dn})}.
\begin{itemize}
\item{} If
\begg \sum_{n=2}^\infty\frac1{D(n)\log n}<\infty, \label{thesum}
\endd then $\{X_n\}$ has finitely many cutpoints almost surely.
\item{} If $D(n)\leq \delta\, n\log n$ $(n\geq n_0)$ for some $\delta>0$
and
$$
\sum_{n=2}^\infty\frac1{D(n)\log n}=\infty,
$$
then $\{X_n\}$ has infinitely many strong cutpoints almost surely.
\end{itemize}
\end{theorem}

\noindent {\bf Remark 2.} Observe that if the sum in (\ref{thesum})
is finite then $\{X_n\}$ has finitely many strong cutpoints. On the
other hand, if the same sum is divergent and $D(n)\leq \delta\,
n\log n$, then $\{X_n\}$ has infinitely many cutpoints as well.

\smallskip
 \noindent {\bf Remark 3.} The condition
$D(n)\leq \delta\, n\log n$ of the second statement is a technical
one, most probably it can be removed.

\smallskip
The above mentioned technical condition
prevent us to establish the following
\begin{conjecture} {\it The number of cutpoints
is finite if and only if the number of strong cutpoints is finite.}
\end{conjecture}

\smallskip
In Section 2 we will present some preliminary results. Sections 3
and 4 are devoted to prove Theorem 1.1. In Section 5 we give some
examples and open problems.

\section{Preliminary results}

\renewcommand{\thesection}{\arabic{section}} \setcounter{equation}{0}
\setcounter{theorem}{0} \setcounter{lemma}{0}

For $0 \leq a\leq b\leq c$ integers define
$$
p(a,b,c)
:=\bp(\min\{j:\ j>m,\ X_j=a\}<\min\{j:\ j>m,\ X_j=c\}\mid X_m=b),
$$
i.e. $p(a,b,c)$ is the probability that a particle
starting from $b$ hits $a$ before $c$.

\bigskip {\bf Lemma A} {\it For} $0 \leq a \leq b \leq c$

$$p(a,b,c)=1-{\frac{D(a,b)}{D(a,c)}}.$$
\noindent {\it Especially, for $n=1,2,\ldots$ we have}

\begin{equation}
p(0,1,n)=1-{\frac{1}{D(0,n)}},\quad p(n,n+1,\infty)=
1-{\frac{1}{D(n)}}.\label{it0}
\end{equation}

It is easy to see that
$$
D(n)=1+D(n+1)U_{n+1},
$$
\begg U_n=\frac{D(n-1)-1}{D(n)},\quad p_n=\frac12\,\,
\frac{D(n)-D(n-1)+1}{D(n)+D(n-1)-1}, \quad n=1,2,\ldots
\label{uni}
\endd

Then observe that for $n\geq m+2$
\begin{eqnarray}
D(m,n)&=&1+\displaystyle{\sum_{j=1}^{n-m-1}\prod_{i=1}^j
U_{m+i}}=D(m)-\sum_{j=n-m}^\infty\prod_{i=1}^j U_{m+i}
=D(m)-U_{m+1}\ldots U_nD(n)\nonumber\\
&=&D(m)-\frac{(D(m)-1)\ldots (D(n-1)-1)}{D(m+1)\ldots D(n-1)}
=D(m)\left(1-\prod_{i=m}^{n-1}\left(1-\frac1{D(i)}\right)\right).
\label{dmn}
\end{eqnarray}

We also define the  number of upcrossings by

\begg \xi(R,n,\uparrow):=\#\{k:\ 0\leq k\leq n,\ X_k=R,\
X_{k+1}=R+1\}.\label{upcr}\endd
 \begg
\xi(R,\uparrow):=\lim_{n\to \infty}\xi(R,n,\uparrow).\endd

It was shown in \cite{CSFR} that

\noindent
{\bf Lemma B} {\it For $R=0,1,2,\ldots $
\begin{equation} {\bf
P}(\xi(R)=L)=\frac{1+2p_R}{2D(R)}
\left(1-\frac{1+2p_R}{2D(R)}\right)^{L-1},\quad L=1,2,\ldots
\label{loc}
\end{equation}
\begin{equation} {\bf
P}(\xi(R,\uparrow)=L)=\frac{1}{D(R)}
\left(1-\frac{1}{D(R)}\right)^{L-1},\quad L=1,2,\ldots
\label{up}
\end{equation}}

It is easy to see that $R$ is a cutpoint if and only if
$\xi(R,\uparrow)=1$. Recall that $R$ is a strong cutpoint if and
only if $\xi(R)=1$.

Denote the (random) set of
\begin{itemize}
\item{} cutpoints by $C$
\item{} strong cutpoints by $C^S$
\end{itemize}

Observe that
\begg
C^S\subseteq C.
\label{sw}
\endd

We present the following exact probabilities.
\begin{lemma} For $k=1,2,\ldots$ we have
\begg \pe( k \in
C)=1-p(k,k+1,\infty)=\frac{1}{D(k)},
\label{pcut}
\endd
\begg \pe(k\in C^S)=\frac{1+2p_k}{2D(k)},
\label{strongcut}
\endd
\begg
\pe(j\in C, k\in C)=\frac{1}{D(j,k+1)D(k)},\quad j<k,
\label{epcut}
\endd
\begg
\pe(j\in C^S, k\in C^S)=\left(\frac12+p_j\right)\left(\frac12+p_k\right)
\frac1{D(j,k)D(k)},\quad j<k.
\label{estrongcut}
\endd
\end{lemma}

\noindent{\bf Proof.} The statements (\ref{pcut}) and (\ref{strongcut})
follow from Lemma B.

To show (\ref{epcut}), we have to observe that after the first arrival
to
$j+1$ the walk has to arrive to $k+1$ without hitting $j$, and from
$k+1$ it must not return to $k$ at all. Formally, by Lemma A
$$
\pe(j\in C, k\in C)=(1-p(j,j+1,k+1))(1-p(k,k+1,\infty))=
\frac{1}{D(j,k+1)D(k)}.
$$

To show (\ref{estrongcut}), observe that after the first hit of $j$ the
walk has to go to $j+1$. Then from $j+1$ it has to hit $k$ before it goes
back to $j$. From $k$ it has to go to $k+1$, and from $k+1$ it must not
return to $k$ at all. Hence again by Lemma A,
$$
\pe(j\in C^S, k\in
C^S)=\left(\frac12+p_j\right)(1-p(j,j+1,k))
\left(\frac12+p_k\right)(1-p(k,k+1,\infty))
$$
$$
=\left(\frac12+p_j\right)\left(\frac12+p_k\right)
\frac{1}{D(j,k)D(k)}.
$$
This completes the proof of the Lemma. $\Box$

\begin{lemma}
For any positive non-decreasing function $G(x), \, x\geq 0$,
the following two sums
$$
\sum_{n=2}^\infty \frac1{G(n)\log n}
\hspace{1.5cm}
\sum_{n=2}^\infty \frac1{G([n\log n])}
$$
are equiconvergent.
\end{lemma}

\noindent
{\bf Proof.} Under the condition of the Lemma we have
$$
\int_n^{n+1}\frac{dx}{G(x)\log x}\leq \frac1{G(n)\log n}\leq
 \int_{n-1}^{n}\frac{dx}{G(x)\log x},
$$
consequently
$$
\int_2^{\infty}\frac{dx}{G(x)\log x}\leq\sum_{n=2}^\infty\frac1{G(n)\log
n}\leq \frac1{G(2)\log 2}+\int_{2}^{\infty}\frac{dx}{G(x)\log x}.
$$

Similarly,
$$
\int_2^\infty\frac{dx}{G(x\log x)}\leq
\sum_{n=2}^\infty \frac1{G([n\log n])}\leq
\frac1{G(1)}+\frac1{G(3)}+\int_2^\infty\frac{dx}{G(x\log x)}.
$$

It remains to show that the integrals
$$
\int^\infty\frac{dx}{G(x)\log x}\quad{\rm and}\quad
\int^\infty\frac{dx}{G(x\log x)}
$$
are equiconvergent. This can be shown by using the substitution $x=y\log
y$ in the first integral above to get
$$
\int^\infty\frac{dx}{G(x)\log x}
=\int^\infty\frac{1+\log y}{\log y+\log\log y}\, \, \frac{dy}{G(y\log
y)}.
$$
Clearly we have
$$
c_1\int^\infty\frac{dy}{G(y\log y)}\leq
\int^\infty\frac{1+\log y}{\log y+\log\log y}\, \, \frac{dy}{G(y\log y)}
\leq c_2\int^\infty\frac{dy}{G(y\log y)}
$$
with some $0<c_1<c_2$, hence the Lemma. $\Box$

\begin{lemma} If $X_n$  and $X_n^*$ are two  NN random walks such
that $X_n^*$ is slower than  $X_n$ then
$$
D(n)\leq D^*(n).
$$
\end{lemma}

\noindent
{\bf Proof.} If $X_n^*$ is slower than  $X_n$ then
$$
U_i=\frac{1/2-p_i}{1/2+p_i}\leq
 \frac{1/2-p^*_i}{1/2+p^*_i}=U_i^*$$
 Hence from the definition  of $D(n)$ we get that
 \begg D(n)\leq D^*(n),   \label{dstar}\endd
proving our Lemma. $\Box$

We will need the following Lemma from Polfeldt \cite{TP}, which is
a particular case of his theorem.

\noindent
{\bf Lemma C} {\it Let $S(x)$ be a slowly varying function such that
$$
\lim_{x\to\infty}\frac{\left(\log S(x)\right)'}{(\log\log x)'}=-\infty
$$
and
$$
\lim_{x\to\infty}\frac{\left(\log S(x)\right)'}{(\log x)'}
\log L(x)=-1
$$
for some normalized differentiable slowly varying function $L(x)$. Then}
$$
\int_x^\infty\frac{S(t)}{t}\, dt\sim S(x)\log L(x), \quad
{\rm as} \quad x\to\infty.
$$

Recall (see Bingham et al. \cite{BGT} page 12-15) that a slowly
varying function $H(x)$ can be represented as
$$
H(x)=a(x)\exp\left(\int_b^x\frac{\varepsilon(t)}{t}\, dt\right),
$$
where $a(x)\to a\neq 0$, $\varepsilon(x)\to 0$ as $x\to\infty$. If
$a(x)=a$, then $H$ is normalized.

Moreover, a differentiable slowly varying function is normalized if and
only if
$$
\frac{xH'(x)}{H(x)}\to 0\quad{\rm as}\quad x\to\infty.
$$

 \setcounter{equation}{0} \setcounter{theorem}{0}
\setcounter{lemma}{0} \setlength{\arraycolsep}{.136889em}

\section{Proof of the convergent part}

We follow the ideas of \cite{JLP}. We have
$\pe(j\in C| k\in C)=1/D(j,k+1)$ for $j<k.$ Observe that it is also the
probability that
\begg \pe(j\in C| k\in C, F_{k+1}) \quad j<k,
\endd
where $F_{k+1}$ is any event determined by the future of the walk
after it reaches $k+1$  for the first time. Let  $C_{j,k}$ be the
set of cutpoints in $(2^j,2^k]$ and $A_{j,k}:=|C_{j,k}|$ the number
of cutpoints in $(2^j,2^k]$. Define

\begg
a_m:=\pe(A_{m,m+1}>0)
\endd
and
\begg
b_m:=\min_{k\in(2^m,2^{m+1}]}
\sum_{i=1}^{2^{m-1}}\frac1{D(k-i,k+1)}
\label{bem}
\endd

On the event that $A_{m,m+1}>0$, let $\ell_m$ be the largest
cutpoint in $C_{m,m+1}.$ We want to give a lower bound for  the
expected number of cutpoints in $(2^{m-1},2^{m+1}]$ by conditioning
on the last cutpoint in $(2^m,2^{m+1}],$ if there is  one:
\begin{eqnarray}
\sum_{j=2^{m-1}+1}^{2^{m+1}} \pe(j\in C) &=&\be(A_{m-1,m+1})\\
\nonumber &\geq & a_m\be(A_{m-1,m+1}|\,A_{m,m+1}>0)\\ \nonumber &=&
a_m\be(\be(A_{m-1,m+1}|\,A_{m,m+1}>0,\ell_m))\\\nonumber &\geq &
a_m\,b_m. \label{ambm}
\end{eqnarray}

It is readily seen that if $p_i\geq 0,\, i=1,2,\ldots$, then
$U_i\leq 1,\, i=1,2,\ldots$ and hence
$$
D(m,n)\leq n-m,
$$
and so
\begg
b_m \geq \sum_{i=1}^{2^{m-1}}\frac{1}{i+1}\geq c\,m \label
{bem2}
\endd
with some $c>0$.

Hence with constants $c$ not necessarily the same on each appearance,
$$
\sum_{m=1}^\infty \bp(A_{m,m+1})=\sum_{m=1}^\infty a_m
\leq \sum_{m=1}^\infty \frac1{b_m}\sum_{j=2^{m-1}+1}^{2^{m+1}}
\bp(j\in C)
$$
$$
\leq \sum_{m=1}^\infty \frac{c}{m}\sum_{j=2^{m-1}+1}^{2^{m+1}}
\frac1{D(j)}
\leq c\sum_{m=1}^\infty\sum_{j=2^{m-1}+1}^{2^{m+1}}
\frac1{D(j)\log j}\leq c\sum_{n=2}^\infty\frac1{D(n)\log n}<\infty,
$$
and by Borel-Cantelli lemma only finitely many of the events
$A_{m,m+1}$ occurs with probability 1, which proves the
convergent part of Theorem 1.1.
$\Box$

 \setcounter{equation}{0} \setcounter{theorem}{0}
\setcounter{lemma}{0} \setlength{\arraycolsep}{.136889em}

\section{Proof of the divergent part}

Let $m_k=[k\log k]$ and
$$
A_k=\{\xi(m_k)=1\}.
$$
We prove that $\bp(A_k\, {\rm i.o.})=1$ which implies the
divergent part of Theorem 1.1. By Lemma 2.1
$$
\bp(A_k)=\frac{1+2p_{m_k}}{2D(m_k)}\geq\frac1{2D([k\log k])},
$$
so by Lemma 2.2,
$$
\sum_k\bp(A_k)=\infty.
$$

For $n>m$ we have
$$
\bp(\xi(m)=1,\, \xi(n)=1)=\left(\frac12+p_m\right)
\frac1{D(m,n)}\left(\frac12+p_n\right)\frac1{D(n)}\leq\frac1{D(m)D(n)H(m,n)}
$$
with
$$
H(m,n)=\frac{D(m,n)}{D(m)}.
$$
It follows from (\ref{dmn}) that
$$
H(m,n)=1-\left(1-\frac1{D(m)}\right)
\ldots\left(1-\frac1{D(n-1)}\right)
$$
$$
\geq 1-\exp\left(-\frac1{D(m)}-\ldots-\frac1{D(n-1)}\right).
$$

Let $\varepsilon>0$ and for given $k$ we split the set $\{\ell>k\}$ into
2 parts. Let
\begg
\ell_1=\min\left\{\ell>k:\, \sum_{i=m_k}^{m_\ell-1}\frac1{D(i)}
\geq\log\frac{1+\varepsilon}{\varepsilon}\right\}.
\label{min}
\endd

\begin{itemize}
\item{(1)}\,\, $\ell\geq \ell_1$,
\item{(2)}\,\, $k<\ell< \ell_1$.
\end{itemize}

In case (1), using that $H(m,n)$ is increasing in $n$ for fixed $m$,
we have for $\ell\geq \ell_1$
$$
\bp(A_k\, A_\ell)=\frac{\bp(A_k)\bp(A_\ell)}{H(m_k,m_\ell)}\leq
\frac{\bp(A_k)\bp(A_\ell)}{H(m_k,m_{\ell_1})}\leq (1+\varepsilon)
\bp(A_k)\bp(A_\ell).
$$

In the case  (2) we use the inequality $1-e^{-u}\geq cu$ for $0\leq
u\leq \log((1+\varepsilon)/\varepsilon)$ with some $c>0$ to get
$$
\bp(A_k\, A_\ell)\leq\frac{\bp(A_k)\bp(A_\ell)}{H(m_k,m_\ell)} \leq
\frac{\bp(A_k)\bp(A_\ell)}{c\,\sum_{i=m_k}^{m_\ell-1}\frac1{D(i)}}
\leq c\,\bp(A_k)\bp(A_\ell) \frac{D(m_\ell)}{m_\ell-m_k}.
$$
Here and in what follows $c,\, c_i$ denote some positive constants, the
values of which might change from line to line.

So we have for $\ell\in$ (2)
$$
\bp(A_k\, A_\ell)\leq \frac{c\,\bp(A_k)}{\ell\log\ell-k\log k},
$$
$$
\sum_{\ell=k+1}^{\ell_1-1}\bp(A_k\, A_\ell)\leq c\,
\bp(A_k)\sum_{\ell=k+1}^{\ell_1-1}\frac1{\ell\log\ell-k\log k}
$$
$$
\leq c\,\bp(A_k)\frac1{\log
k}\sum_{\ell=k+1}^{\ell_1-1}\frac1{\ell-k}\leq
c\,\bp(A_k)\frac{\log\ell_1}{\log k}.
$$
Now we show that
\begg
\frac{\log\ell_1}{\log k}\leq \gamma
\label{square}
\endd
with some positive constant $\gamma$ depending only on $\varepsilon$. We
know from (\ref{min}) that for $\ell\in (2)$ we have
\begin{equation}
\sum_{i={m_k}}^{m_{\ell}-1}\frac1{D(i)}<\log\frac{1+\varepsilon}{\varepsilon}.
\label{log}
\end{equation}
We show that this implies that for large $k$ we have $\ell<k^{\gamma}$
with $\gamma>(1+\varepsilon/\varepsilon)^\delta$.
If we assume the contrary that $\ell\geq k^{\gamma}$, then
$$
\sum_{i=m_k}^{m_{\ell}-1}\frac1{D(i)}\geq
\frac1{\delta}\sum_{i=m_k}^{m_{\ell}-1}\frac1{i\log i}
\sim\frac1{\delta}(\log\log(m_\ell-1)-\log\log m_k)
$$
$$
\sim\frac1{\delta}\log\frac{\log(\ell\log\ell)}{\log(k\log k)}
\geq\frac1{\delta}\log\gamma>\log \frac{1+\varepsilon}{\varepsilon}
$$
which contradicts to (\ref{log}). Hence $\ell_1-1\leq k^{\gamma}$,
implying (\ref{square}).

Consequently,
$$
\sum_{\ell\in (2)}\bp(A_k\, A_\ell)\leq c\bp(A_k)
$$

Assembling these estimations, we have
$$
\sum_{k=1}^N\sum_{\ell=k+1}^N\bp(A_k\, A_\ell)\leq
(1+\varepsilon)\sum_{k=1}^N\sum_{\ell=k+1}^N\bp(A_k)\bp(A_\ell)
+c\sum_{k=1}^N\bp(A_k)
$$
Since $\varepsilon>0$ is arbitrary, Borel-Cantelli lemma implies
$\bp(A_k\, {\rm i.o.})=1$. $\Box$

\setcounter{equation}{0} \setcounter{theorem}{0}
\setcounter{lemma}{0} \setlength{\arraycolsep}{.136889em}

\section{Conclusions and open problems}

Our results are formulated in terms of the sequence $\{D(\cdot)\}$
but it would be much more natural to formulate them in terms of the
sequence $\{p_i\}$. Even though we have an explicit expression of
$D(\cdot)$ in terms of $\{p_i\}$, usually it is not easy to see the
asymptotics of $D(\cdot)$ and whether the sum in Theorem 1.1 is
convergent or divergent by looking at $\{p_i\}$ only. Therefore we want
to give some examples.

In \cite{CSFR} we have shown

\noindent
{\bf Example 1.} {\it If $p_k=B/4k$ with $B>1$, then}
$$
D(i)\sim \frac{i}{B-1}
$$
{\it as $i\to\infty$. Consequently, by {\rm Theorem 1.1} there are
infinitely many strong cutpoints a.s.}

\medskip
It was shown also in \cite{CSFR}

\noindent
{\bf Example 2.} {\it If $p_k=\Lambda(K,k,B)/4$ with $K\geq 2$ and
$B>1$, then}
$$
D(i)\sim \frac{i\log i\log\log i\ldots \log_{K-1}i}{B-1}
$$
{\it as $i\to\infty$. Consequently, by {\rm Theorem 1.1} we have
finitely many cutpoints a.s.}

Recall that the case $K=2$ corresponds to Theorem D.

Of course, if the NN walk is quicker than the walk in Example 1,
(e.g. $p_k=c\,k^{-\alpha}$ with $\alpha<1$), then we have infinitely
many strong cutpoints a.s. On the other hand, if the NN walk is
slower than the walk in Example 2, then we have finitely many
cutpoints a.s.

The above two examples show that the jump from finitely many to
infinitely many cutpoints is for
$$
p_k=\frac1{4}\left(\frac1k+\frac1{kf(k)}\right)
$$
with some $f(k)\to\infty$.
It is not hard to show that if $f(k)=(\log k)^{\alpha}$ with
$0<\alpha<1$, then we still have finitely many cutpoints a.s. Now we
show a more precise result which implies this one.

\begin{theorem}
Let $\{X_n\}$ be an NN random walk with
$$
p_k=\frac{1}{4}\left(\frac1k+\frac{1}{k(\log\log k)^{\beta}}\right),
$$
then we have finitely many cutpoints a.s. if $\beta>1$ and infinitely
many strong cutpoints a.s. if $\beta\leq 1$.
\end{theorem}

\noindent
{\bf Proof.}
Let
$$
r_k=\prod_{i=1}^k U_i,\quad t_k=\sum_{i=k}^\infty r_i,\qquad
k=1,2,\ldots
$$
Then it is easy to see that $D(k)=t_k/r_k$.

For $k\to\infty$ we obtain
$$
U_k=\frac{1-2p_k}{1+2p_k}=\exp(-4p_k+O(p_k^2))
$$
and
$$
r_k=\exp(-4\sum_{i=1}^k p_i+O(\sum_{i=1}^k p_i^2))
=\exp\left(-\sum_{i=3}^k\left(\frac 1{i}+\frac1{i(\log\log
i)^\beta}\right)+O(1)\right)
$$
$$
=\exp\left(-\int_{3}^k\left(\frac1{u}+\frac1{u(\log\log
u)^\beta}\right)\, du+O(1)\right).
$$
Hence
\begg
\frac{c_1}{k}\exp\left(-\int_{3}^k\frac{du}{u(\log\log
u)^\beta}\right)\leq r_k\leq
\frac{c_2}{k}\exp\left(-\int_{3}^k\frac{du}{u(\log\log
u)^\beta}\right)
\label{rk}
\endd
with some positive constants $c_1,\, c_2$. Consequently,
\begg
c_1\sum_{j=k}^\infty\frac1{j}\exp\left(-\int_{3}^j\frac{du}{u(\log\log
u)^\beta}\right)\leq t_k\leq
c_2\sum_{j=k}^\infty\frac1{j}\exp\left(-\int_{3}^j\frac{du}{u(\log\log
u)^\beta}\right).
\label{tk}
\endd
Moreover,
$$
\sum_{j=k}^\infty\frac1{j}\exp\left(-\int_{3}^j\frac{du}{u(\log\log
u)^\beta}\right)=
\int_k^\infty\frac1{y}\exp\left(-\int_{3}^y\frac{du}{u(\log\log
u)^\beta}\right)+O(1).
$$
To find the asymptotics of the above integral, we will apply Lemma C,
with
$$
S(y)=\exp\left(-\int_{3}^y\frac{du}{u(\log\log
u)^\beta}\right)
$$
and
$$
L(x)=e^{(\log\log x)^\beta}.
$$
Choosing $S(\cdot)$ and $L(\cdot)$ as above, all the conditions of Lemma
C are met and we conclude that
$$
\int_k^\infty\frac1{y}\exp\left(-\int_{3}^y\frac{du}{u(\log\log
u)^\beta}\right)\sim (\log\log k)^\beta
\exp\left(-\int_{3}^k\frac{du}{u(\log\log
u)^\beta}\right).
$$
From (\ref{rk}) and (\ref{tk}) we obtain
$$
c_1\, k\, (\log\log k)^\beta\leq D(k)=\frac{t_k}{r_k}\leq
c_2\, k\, (\log\log k)^\beta.
$$
This combined with Theorem 1.1 proves Theorem 5.1. $\Box$

Lemma 2.4 easily implies the following
\begin{corollary}
If $\{X_n\}$ is an NN random walk with
$$
p_k\leq \frac14\left(\frac1{k}+\frac{1}{k(\log\log k)^\beta}\right)
$$
and $\beta>1$, then $X_n$ has finitely many cutpoints a.s.

On the other hand, if $\{X_n\}$ is an NN random walk with
$$
p_k\geq \frac14\left(\frac1{k}+\frac{1}{k(\log\log k)^\beta}\right)
$$
and $\beta\leq 1$, then $X_n$ has infinitely many strong cutpoints a.s.
\end{corollary}

Now we present some related open problems.
\begin{itemize}
\item{(1)} It would be interesting to know whether Theorem 1.1 also
holds for the number of sites with $\xi(R)=a$ or $\xi(R,\uparrow)=a$ for
any fixed integer $a>1$, i.e. whether we have the same criteria for
$\{\xi(R)=a\}$ and $\{\xi(R,\uparrow)=a\}$ to occur infinitely often
almost surely for any positive integer $a$.
\item{(2)} Call the site $R$ a {\it weak cutpoint}\, if for some $k$, we
have $X_k=R$, $X_i\leq R,\, i=0,1,\ldots, k-1$ and $X_i\geq R,\,
i=k+1,k+2,\ldots$. One would like to know whether Theorem 1.1 can be
extended for the number of weak cutpoints.
\item{(3)} It would be interesting to know whether Theorem 1.1 holds for
cutpoints with a given local time, i.e. for
$\{\xi(R)=a,\,\xi(R,\uparrow)=1\}$, or in general
$\{\xi(R)=a,\,\xi(R,\uparrow)=b\}$ infinitely often almost surely, with
positive integers $a,b$.
\item{(4)} Theorem B gives limsup behavior of the local time. One might
ask how does it change if we want to consider the limsup of the local
time restricted to the cutpoints.
\end{itemize}

\end{document}